\documentclass[12pt]{amsart}
\usepackage{times}
\usepackage{amsfonts}
\usepackage{amsthm}
\usepackage{amsmath}
\advance \topmargin by -\headheight
\advance \topmargin by -\headsep
\evensidemargin \oddsidemargin
\newtheorem{theorem}{Theorem}[section]
\newtheorem{lemma}[theorem]{Lemma}

\newtheorem{remark}[theorem]{Remark}
\newtheorem{example}[theorem]{Example}

\begin{document}

\title{Fractal Entropies and Dimensions for Microstate Spaces, II}

\author{Kenley Jung}

\address{Department of Mathematics, University of California,
Berkeley, CA 94720-3840,USA}

\email{factor@math.berkeley.edu}
\thanks{Research supported by the NSF Graduate Fellowship Program}

\begin{abstract}For a selfadjoint element $x$ in a tracial von Neumann
algebra and $\alpha = \delta_0(x)$ we compute bounds for $\mathbb
H^{\alpha}(x),$ where $\mathbb H^{\alpha}(x)$ is the free Hausdorff
$\alpha$-entropy of $x.$ The bounds are in terms of $\int \int_{\mathbb
R^2 -D} \log |y-z| \, d\mu(y) d\mu(z)$ where $\mu$ is the Borel measure on
the spectrum of $x$ induced by the trace and $D \subset \mathbb R^2$ is
the diagonal.  We compute similar bounds for the free Hausdorff entropy of
a free family of selfadoints.  \end{abstract} \maketitle

\section*{introduction}

[1] introduced fractal geometric entropies and dimensions for Voiculescu's
microstate spaces ([3], [4]).  One can associate to a finite set of
selfadjoint elements $X$ in a tracial von Neumann algebra and an $\alpha
>0$ an extended real number $\mathbb H^{\alpha}(X) \in [-\infty, \infty].$
$\mathbb H^{\alpha}(X)$ is a kind of asymptotic logarithmic
$\alpha$-Hausdorff measure of the microstate spaces of $X.$ One can also
define a free Hausdorff dimension of $X,$ denoted by $\mathbb H(X),$ which
is related to $\mathbb H^{\alpha}(X)$ in the same way that Hausdorff
dimension is related to the critical value of Hausdorff measures. $\mathbb
H^{\alpha}$ can be regarded as an interpolated version of Voiculescu's
free entropy $\chi$ in the sense that if $X$ consists of $n$ selfadjoints,
then $\mathbb H^n(X) = \chi(X) + \frac{n}{2} \log ( \frac{2n}{\pi e}).$
This connection seems perfectly natural since $\chi$ is defined in terms
of Lebesgue measure and Hausdorff $n$ measure is a normalization of
$n$-dimensional Lebesgue measure.

In [3] Voiculescu establishes an equation for $\chi(x)$ where $x$ is a
selfadjoint operator.  He shows that if $\mu$ is the Borel measure on
$sp(x)$ induced by the tracial state, then the free entropy of $x$ is a
normalization of the logarithmic energy of $\mu,$ i.e.,

\[ \chi(x) = \int \int \log |y-z| \, d\mu(y) d\mu(z) + \frac{3}{4} + 
\frac{1}{2} \log 2 \pi.\]

\noindent Moreover, Voiculescu showed in the same work that if 
$\{x_1,\ldots, x_n\}$ is a free family of selfadjoints, then 
$\chi(x_1,\ldots, x_n) = \chi(x_1) + \cdots + \chi(x_n).$ 

It is natural to wonder whether similar properties hold for the free
Hausdorff $\alpha$-entropy.  The strongest statement in this direction
might go as follows.  If $x$ is selfadjoint and $\alpha$ is the free
Hausdorff dimension of $x,$ then

\[ \mathbb H^{\alpha}(x) = \int \int_{\mathbb R^2-D} \log |y-z| \, 
d\mu(y) d\mu(z) +K_{\alpha},\]  

\noindent where $K_{\alpha}$ is some constant dependent on $\alpha$ and
$D$ is the diagonal line in $\mathbb R^2.$ Using Voiculescu's
strengthened asymptotic freeness results, free additivity would follow.  
If this is too much to ask for, then one might hope to show that $\mathbb
H^{\alpha}(x)$ is bounded in terms of $\int \int_{\mathbb R^2-D} \log
|y-z| \, d\mu(y) d\mu(z).$ Presumably, this would be followed by showing
that a free family of selfadjoints has Hausdorff entropy proportional to
the sums of the free Hausdorff entropies of each element.  
Unfortunately, due to technical difficulties and some fundamental
differences between $\mathbb H^{\alpha}$ and $\chi,$ neither equations
nor estimates of these kinds were present in [1].

The present work is an addendum to [1] where we fill in this gap
by showing the weaker of the two proposed problems, namely for
$\alpha = \mathbb H(x)$ there exist constants $K_1, K_2$ dependent only on
$\alpha$ such that

\[ K_1 \leq \int \int_{\mathbb R^2 -D} \log |y-z| \, d\mu(y) d\mu(z) -
\mathbb H^{\alpha}(x) \leq K_2.\]

\noindent Moreover, we show that these bounds promote to ones for a freely
independent, finite family of selfadjoints.

The techniques are very much in the spirit of those in part 1.  Because
all the microstate spaces are naturally associated to locally isometric
spaces all the Hausdorff quantities can be bound with strong packing
estimates and these in turn can be computed by results of Mehta.  The only 
new aspects involve sharpening the aforementioned methods to arrive at 
tighter estimates.  
	
There are five short sections.  The first is a list of notation.  The
second is a brief collection of properties we will use about microstates
for a single selfadjoint.  The third and fourth sections are the upper and
lower bounds, respectively, on the free Hausdorff entropy of a single
selfadjoint, and the fifth section deals with the free situation.

\section{Notation} Throughout suppose $M$ is a von Neumann algebra with
faithful, tracial state $\varphi.$ Suppose $x = x^* \in M$ and $\mu$ is
the Borel measure on $sp(x)$ induced by $\varphi.$ $\alpha = \mathbb
H(x),$ which by [1] and [4] is just $\delta_0(x) = \delta(x).$ Suppose $R
> \|x\|.$ We maintain the standard notation introduced in [4] for the
microstate spaces.  $M^{sa}_k(\mathbb C)$ denotes the set of selfadjoint
$k\times k$ complex matrices and $(M^{sa}_k(\mathbb C))^n$ is the space of
$n$-tuples with entries in $M^{sa}_k(\mathbb C).$ All metric quantities
for the microstate spaces will be taken with respect to the norm $| 
\cdot|_2$ on $(M^{sa}_k(\mathbb C))^n$ given by $|(a_1,\ldots,a_n)|_2 =
(\sum_{j=1}^n tr_k(a_j^2))^{\frac{1}{2}}$ where $tr_k$ is the tracial
state on $M_k(\mathbb C).$ vol denotes Lebesgue measure on
$(M^{sa}_k(\mathbb C))^n$ with respect to the Hilbert space norm $\sqrt{k}
\cdot | \cdot |_2$ and $L_k$ denotes the Lebesgue measure of the ball of
radius $\sqrt{k}$ in $\mathbb R^{k^2}.$ $D$ will denote the diagonal line
in $\mathbb R^2.$

\section{Microstates for a single selfadjoint}

Because the von Neumann algebra generated by $x$ is hyperfinite, the
microstate space for a single selfadjoint $x$ is obtained by taking
unitary orbit of one well-approximating microstate.  Since the estimates
involve the entropy (and not dimension), we will need a sharper handle on
such `well-approximating' microstates and define two kinds of microstates
for $x$: $A_k$ and $B_k.$ The $A_k$ will be used for the upper
bound and the $B_k$ will be used for the lower bound.

Write $\mu = \sigma + \nu$ where $\sigma$ is the atomic part of $\mu$ and
$\nu$ is the diffuse part of $\mu.$ $\sigma = \sum_{i=1}^s c_i
\delta_{r_i}$ for some $s \in \mathbb N \cup \{0\} \cup \{\infty\},$ $c_i
\geq c_{i+1} >0,$ and where for $i \neq j,$ $r_i \neq r_j.$ Set $c =
\nu([a,b])$ where $sp(x) \subset [a,b].$ Because $\nu$ is diffuse for each
$k$ and $1 \leq j \leq [ck]$ there exists a largest number $\lambda_{jk}
\in [a,b]$ satisfying $\nu([a, \lambda_{jk}]) = \frac{j}{k}.$ 

\subsection{$A_k$ microstates for the upper bound}

For each $k$ denote by $A_k$ the $k \times k$ diagonal matrix obtained by
filling in the the first $[ck]$ entries with $\lambda_{1k}, \ldots,
\lambda_{[ck]k}$ and the last $\sum_{j=1}^s [c_jk]$ diagonal entries
filled with $r_1$ repeated $[c_1k]$ times, $r_2$ repeated $[c_2k]$ times,
etc., in that order.  Fill in the remaining $k - [ck] - \sum_{i=1}^s
[c_ik]$ terms with $0$'s.  Observe that for any $m \in \mathbb N$ and
$\gamma > 0$ $A_k \in \Gamma_R(x;m,k,\gamma)$ for $k$ sufficiently large.
Also observe that if $\epsilon >0,$ and $a_{1k}, \ldots, a_{kk}$ are the
eigenvalues of $A_k$ ordered from least to greatest and according to
multiplicity, then

\[\lim_{k \rightarrow \infty} k^{-2} \log \Pi_{1 \leq i < j \leq k}
((a_{ik} - a_{jk})^2 + \epsilon) = \int \int \log (|y-z|^2 + \epsilon) \,
d\mu(y) d\mu(z). \]

\noindent This follows by writing each term in the limit as the integral
of an obvious simple function $f_k(y,z)$ defined on $[a,b]^2.$ $\langle
f_k \rangle_{k=1}^{\infty}$ will be a sequence uniformly bounded by
$\max\{ |\log \epsilon|, |\log (b-a)|\}$ and $f_k(s,t) \rightarrow \log
((y-z)^2 + \epsilon $ a.e. $\mu \times \mu.$ It then follows from
Lebesgue's Dominated Convergence Theorem that

\[k^{-2} \log \Pi_{1 \leq i < j \leq k} ((a_{ik} - a_{jk})^2 + \epsilon) =
\int \int f_k(y,z) \, d\mu(y) d\mu(z)  \rightarrow \int \int \log (|y-z|^2
+ \epsilon) \, d\mu(y) d\mu(z). \]

\subsection{$B_k$ microstates for the lower bound}

The $B_k$ are defined only when $s \geq 1,$ i.e., $\mu$ has a nontrivial
atomic part, and when $\int \int_{\mathbb R^2 -D} \log |y-z| \, d\mu(y)
d\mu(z)$ is finite.  In this case $B_k$ will be the $k\times k$ diagonal
matrix obtained by first adding $r_1$ $[c_1 k] - \sqrt{k}$ times, followed
by adding $r_2$ $[c_2 k]$ times, then $r_3$ added $[c_3k]$ times, and
continuing in this way. This process will terminate for there exists a
maximum value $N_k \in \mathbb N$ dependent on $k$ for which if $j> N_k$
then $[c_jk] =0.$ Now recall the $\lambda_{jk}$ defined in the preceding
paragraph.  For each $1 \leq m \leq N_k$ find the largest $\lambda_{jk}$
less than or equal to $r_j$ and the smallest $\lambda_{jk}$ greater than
or equal to $r_j.$ Denote by $R_k$ the set of all such $\lambda_{jk}.$
Observe that $\#R_k \leq 2 N_k.$ Fill in the remaining entries of $B_k$
with $\{\lambda_{2k},\ldots, \lambda_{([ck]-1)k}\} - R_k,$ ordered from
greatest to least.  This leaves a remaining $F_k$ entries to fill in where
$F_k \leq 2 N_k + \sum_{j= N_k+1} c_j k.$ Fill these entries with $B +3 +
\frac{1}{F_k}, B +3 + \frac{2}{F_k}, \ldots, B + 4.$ For any given $m \in
\mathbb N$ and $\gamma >0$, $B_k \in \Gamma_{R+4}(x;m,k,\gamma)$ for $k$
sufficiently large.  Let $b_{1k}, \ldots, b_{kk}$ be the eigenvalues of
$B_k$ ordered from least to greatest and with respect to multiplicity.

Denote by $S_k$ the set of all $2$-tuples $(i,j)$ such that $1 \leq i < j 
\leq k,$ and $b_{ik} = b_{jk}$; this can only happen when $b_{ik} = 
b_{jk}$ is one of the atoms $r_1,\ldots, r_{N_k}$ of $\sigma.$  Denote by 
$W_k$ the set of all $2$-tuples $(i,j)$ such that $1 \leq i <j \leq k,$ 
and $(i,j) \notin S_k.$  I claim that 

\begin{eqnarray} \liminf_{k \rightarrow \infty} k^{-2} \cdot 
\sum_{(i,j) \in W_k} \log (b_{ik} - b_{jk})^2 \geq \int \int_{\mathbb R^2 
-D} \log |y-z|\, 
d\mu(y)d\mu(z).\end{eqnarray}

\noindent Fix $k.$ Define $X_k$ to consist of all $(i,j) \in W_k$ such
that $1 \leq i < j \leq k,$ either $b_{ik}$ or $b_{jk}$ does not belong
to $\{r_1,\ldots, r_{N_k}\},$ and satisfying the condition that $a \leq
b_{ik}, b_{jk} \leq b.$ If $b_{ik}$ is not an element of $\{\lambda_{1k},
\ldots, \lambda_{[ck]k}\},$ then denote by $R(i)$ the smallest element in
$\{\lambda_{1k}, \ldots, \lambda_{[ck]k}\}$ larger than $b_{ik}$ and
$L(i)$ the largest element in $\{\lambda_{1k}, \ldots, \lambda_{[ck]k}\}$
smaller than $b_{ik}.$ If neither $a_{ik}$ nor $b_{jk}$ belong to
$\{r_1,\ldots, r_{N_k}\}$ then set $f_{ij}$ to be $\log |b_{ik} -
b_{jk}|^2$ times the characteristic function over the set $ ((b_{ik},
R(i)) - \{r_j: 1 \leq j \leq s\})  \times ((L(j),b_{jk}) - \{r_j: 1 \leq
j \leq s\}).$ Otherwise either $b_{ik} \in \{r_1,\ldots, r_{N_k}\}$ or
$b_{jk} \in \{r_1,\ldots,r_{N_k}\}$ but not both.  In the former case
define $f_{ij}$ to be $\frac{1}{c_p k} \cdot \log |a_{jk} - r_p|^2$
times the characteristic function over $\{r_p\} \times \left ((L(j),
b_{jk}) - \{r_j: 1 \leq j \leq s\} \right)$ where $r_p = b_{ik}.$ Observe
that $r_p < L(j).$ In the latter case define $f_{ij}$ to be
$\frac{1}{c_p k} \cdot \log |b_{ik} - r_p|^2$ times the characteristic
function over $\left ((b_{ik}, R(i)) - \{r_j: 1\leq j \leq s\} \right)
\times \{r_p\}$ where $r_p = b_{jk}.$ Observe that $R(i) < r_p.$ Finally,
set

\[ g_k(y,z) = \sum_{1 \leq i < j \leq N_k} \log |r_i - r_j| \cdot
\chi_{(\{r_i\} \times \{r_j\})} (y,z) + \sum_{(i,j) \in X_k} f_{ij}(y,z).  
\]

\noindent Because for each $(i,j) \in X_k,$ $\int \int f_{ij}(y,z) \, 
d\mu(y) d\mu(z) = k^{-2} \cdot \log (b_{ik} - b_{jk})^2$ it is clear that 

\begin{eqnarray} \int \int_{r<t} g_k(y,z) \, d\mu(y) d\mu(z) \leq \sum_{1
\leq i < j \leq k} k^{-2} \log (b_{ik} - b_{jk})^2 - \sum_{1 \leq i
< j \leq N_k} k^{-2} \log \left (\frac {j-i}{k} \right).
\end{eqnarray}

\noindent Notice that the second term in the sum on the right hand side
above converges to $0$ as $k \rightarrow 0$ since $\lim_{k \rightarrow
\infty} \frac{N_k}{k} = 0.$  $|g_k(y,z)| \leq \max \{ |\log
(y-z)^2|, |\log (b-a)|\}$ for any $k.$  

Moreover, $g_k(y,z) \rightarrow \log (y-z)^2$ a.e.  $\mu \times \mu$ on
the region $E$ strictly above $D.$ To see this suppose $\epsilon,L >0$ and
choose $n \in \mathbb N$ so large that $\sum_{j=n}^s c_j < \epsilon.$
Write $K = \text{support} \hspace{.03in} (\nu) - \{r_j : 1\leq j \leq
s\}.$ There exists a $\delta >0$ such that the $\mu \times \mu$ measure of
$\{(y,z) \in \mathbb R^2:  0 < |y-z| < \delta\}$ is no greater than
$\epsilon.$ Define $T$ to be the intersection of the region strictly above
the line $y = z + \delta$ and

\[  \left (\{r_j: 1\leq j \leq n\} \cup K \right)^2 .\]

\noindent $(\nu \times \nu)(E -T) < \epsilon.$ Define for each $k,$ $G_k =
\{ 1 \leq i \leq [ck] : |\lambda_{ik} - \lambda_{jk}| < \frac{\delta}{L},
j = i \pm 1\}.$ It is each to show that $\lim_{k \rightarrow \infty}
\#G_k/ck =1.$ If $H_k$ consists of all $b_{jk}$ such that $b_{jk} =
\lambda_{ik}$ for some $i \in G_k,$ then $\lim_{k \rightarrow \infty}
\#H_k / ck = 1.$ Write $I_k$ for the union of all Cartesian products of
the form $((b_{ik},R(i)) - \{r_j: 1\leq j \leq s\}) \times ((L(j), b_{jk})
- \{r_j : 1 \leq j \leq s\})$ , $\{r_p\} \times ((L(j), b_{jk}) - \{r_j :
1\leq j \leq s\})$ ( when $r_p = b_{ik}$), and $((b_{ik},R(i)) - \{r_j: 1
\leq j \leq s\}) \times \{r_p\}$ (when $r_p = b_{ik}$) where $(i,j) \in
X_k \cap H_k^2.$ Take all the union of all these sets with $\{r_j: 1 \leq
j \leq n\}^2$ and call the resultant set $I_k.$ Observe that $\lim_{k
\rightarrow \infty} (\mu \times \mu)(T - I_k) = 0$ and that on $I_k$
$g_k(y,z)$ and $\log(y-z)^2$ differ by no more than $\frac{\delta}{L}
\cdot \max \{ |\log \delta|, |\log(b-a)|\}.$ As $L$ was arbitrary, it
follows that $g_k(y,z)$ converges to $\log(y-z)^2$ almost everywhere $\mu
\times \mu$ on $T.$ Since $\epsilon >0$ was arbitrary and $(\mu \times
\mu)(E-T) < \epsilon,$ $g_k(y,z)$ converges to $\log (y-z)^2$ on $E$ a.e.  
$\mu \times \mu.$

It now follows from Lebesgue's Dominated Convergence Theorem that

\[ \int \int_{E} g_k(y,z) \, d\mu(y) d\mu(z) \rightarrow \int \int_{E} 
\log (y-z)^2 \, d\mu(y) d\mu(z) =
\int \int_{\mathbb R^2 -D} \log |y-z| \, d\mu(y)d\mu(z).\]

\noindent Combining this with the preceding inequality gives (1).

\section{upper bound}

\begin{lemma} $\mathbb H^{\alpha}(x) \leq \int \int_{\mathbb R^2 - D} \log |s-t| d\mu(s)d\mu(t) + \log 16 +
\frac{1}{4}.$ \end{lemma}

\begin{proof} Suppose $\epsilon, t >0$ are given.  There exist $m \in \mathbb N$ and $\gamma >0$ such that for any
$k \in \mathbb N$ and $A, B \in \Gamma(x;m,k,\gamma)$ there exists a unitary $u$ satisfying $|uAu^* -B |_2 < t.$
Consider the sequence $\langle A_k \rangle_{k=1}^{\infty}$ constructed in Lemma 2.1.  For sufficiently large $k,$
$A_k \in \Gamma(x;m,k,\gamma)$ and moreover,

\begin{eqnarray*} H^{\alpha k^2}_{\epsilon}(\Gamma(x;m,k,\gamma)) \leq H^{\alpha k^2}_{\epsilon}(\Theta_t(A_k))
\leq K_{\frac{\epsilon}{2}}(\Theta_t(A_k)) \cdot \epsilon^{\alpha k^2} & \leq & P_{\frac{\epsilon}{4}}
(\Theta_t(A_k)) \cdot \epsilon^{\alpha k^2} \\ & \leq & \text{vol} \hspace{.03in} (\Theta_{t+ \epsilon/4 }(A_k))
\cdot \frac{4^{k^2} \epsilon^{(\alpha -1)k^2}}{L_k}.  \end{eqnarray*}

\noindent By Lemma 4.2 of [3] $\text{vol} \hspace{.03in} (\Theta_{t + \epsilon/4}(A_k))$ is dominated by

\[ k^{k/2} \cdot \epsilon^k \cdot \Gamma\left (\frac{k}{2} + 1\right)^{-1} \cdot (1 + 2\alpha)^{k(k-1)/2} \cdot
e^{2k^2 \epsilon} \cdot \pi^{k^2/2} \cdot 2^{k (k-1)/2} \cdot (\Pi_{j=1}^k j!)^{-1} \cdot \Pi_{1\leq i <j \leq k}
((a_{ik} - a_{jk})^2 + \epsilon).\]

\noindent Here $\alpha \in (0, \frac{1}{2})$ is the unique number such that $(\alpha+ 2 \alpha^2)(\alpha
+2)^{-1})^{1/2} = \frac{t}{\epsilon} + \frac{1}{4}.$ Thus, using Lemma 4.4 of [3] and the preliminary remarks in
Section 2 it follows that

\begin{eqnarray*} \mathbb H^{\alpha}_{\epsilon}(x) & \leq & \mathbb H^{\alpha}_{\epsilon}(x;m,\gamma) \\ & \leq &
\limsup_{k \rightarrow \infty} k^{-2} \cdot \left [\log (\text{vol} \hspace{.03in}(\Theta_{t + \epsilon/4}(A_k)) -
\log L_k \right] + \log 4 + (\alpha -1) \log \epsilon \\ & \leq & \limsup_{k \rightarrow \infty} k^{-2} \cdot \log
\Pi_{1 \leq i < j \leq k} ((a_{ik} - a_{jk})^2 + \epsilon) + \log 16 + (\alpha -1) \log \epsilon + 2 \epsilon +
\frac{1}{4} \\ & \leq & \frac{1}{2} \cdot \int \int_{\mathbb R^2} 
\log(|s-t|^2 + \epsilon) \, d\mu(s) d\mu(t)+ (\alpha -1) \log \epsilon +
\log 16 + \frac{1}{4} + 2 \epsilon \\ & \leq &\frac{1}{2} \cdot \int 
\int_{\mathbb R^2 -D} \log(|s-t|^2 + \epsilon) \, d\mu(s)
d\mu(t) + \log \epsilon \cdot (\mu \times \mu)(D)+ (\alpha -1) \log 
\epsilon \\ & & + \log 16 + \frac{1}{4} + 2 \epsilon \\ &
\leq & \frac{1}{2} \cdot \int \int_{\mathbb R^2 -D} \log(|s-t|^2 + \epsilon) \, d\mu(s)d\mu(t) + \log 16 +
\frac{1}{4} + 2 \epsilon\\ \end{eqnarray*}

\noindent Forcing $\epsilon \rightarrow 0$ we arrive at the desired conclusion.  
\end{proof}

\section{lower bound}

\begin{lemma} $\lim_{k \rightarrow \infty}
k^{-2} \cdot \log \left [\Pi_{j=1}^k \frac{\Gamma(j+1)
\Gamma(j)^2}{\Gamma(k+j)} \right] = - \log 4.$ \end{lemma}

\begin{proof} Using Lemma 4.4 in [3],
\begin{eqnarray*} \lim_{k \rightarrow \infty} k^{-2} \cdot \log \left
[\Pi_{j=1}^k \frac{\Gamma(j+1) \Gamma(j)^2}{\Gamma(k+j)} \right] & = &
\lim_{k \rightarrow \infty} k^{-2} \cdot \log \left [\frac {\Pi_{j=1}^k
\Gamma(j)^4}{\Pi_{j=1}^{2k} \Gamma(j)} \right] \\ & = & 4 \cdot \lim_{k
\rightarrow \infty} \left ( k^{-2} \cdot \log \Pi_{j=1}^k \Gamma(j) -
2^{-1} \log k \right) + \\ & & 4 \cdot \lim_{k\rightarrow \infty} \left (
-(2k)^{-2} \cdot \log \Pi_{j=1}^{2k} \Gamma(j) + 2^{-1} \log 2k \right)-
\\ & & 2 \log 2 \\ & = & -3 + 3 - 2 \log 2 \\ & = & - \log 4.
\end{eqnarray*} \end{proof} 

\begin{lemma} $\mathbb H^{\alpha}(x) \geq - \delta_0(x) \log 2 -
\frac{1}{2} \log 288e + \frac{3}{4} + \int \int_{\mathbb R^2 -D} \log |s-t| \,
d\mu(s) d\mu(t).$ \end{lemma} 

\begin{proof} Note that the inequality trivially holds when $\mathbb H(x)  =1,$ i.e., when $x$ has no eigenvalues.  
This follows from Proposition 4.4 in [3] and Lemma 3.7 of [1].  Also observe that the desired
inequality is vacuously satisfied when the integral in question is $- \infty$ Thus, we assume without loss of
generality, that $x$ has a nontrivial point spectrum and $\int \int_{\mathbb R^2 -D} \log |y-z| \, d\mu(y) d\mu(z)
> -\infty.$ 

Denote by $G$ the group of diagonal unitaries and $\mathbb R^k_{<}$ to be
the set of all $(t_1,\ldots,t_k) \in \mathbb R^k$ such that $t_1 < \cdots
< t_k.$ There exists a map $ \Phi : M^{sa}_k(\mathbb C) \rightarrow U_k/G
\times \mathbb R^k_{<}$ defined almost everywhere on $M^{sa}_k(\mathbb C)$
such that for each $x \in M^{sa}_k(\mathbb C)$ $\Phi(x) = (h, z)$ where
$z$ is a diagonal matrix with real entries satisfying $z_{11} < \cdots <
z_{kk}$ and $h$ is the image of any unitary $u$ in $U_k/G$ satisfying
$uzu^* = x.$ By results of Mehta ([2]) the map $\Phi$ induces a measure
$\mu$ on $U_k/G \times \mathbb R^k_<$ given by $\mu(E) = \text{vol}
(\Phi^{-1}(E))$ and moreover, \[ \mu = \nu \times D_k \cdot \int_{\mathbb
R^k_{<}} \Pi_{i < j} (t_i - t_j)^2 \, dt_1 \cdots dt_k \] \noindent where
$D_k =\frac{\pi^{k(k-1)/2}}{\Pi_{j=1}^k j!}$ and $\nu$ is the normalized
measure on $U_k/G$ induced by Haar measure on $U_k.$ Write
$\Theta_{\epsilon}(B_k)$ for the $|\cdot|_2$ $\epsilon$-neighborhood of
the unitary orbit of $B_k$ (as defined in subsection 2.2) and 
$\Theta(B_k)$ for the unitary orbit of
$B_k.$ A matrix will be in $\Theta_{\epsilon}(B_k)$ iff the sequence
obtained by listing its eigenvalues in increasing order and according to
multiplicity, differs from the similar sequence obtained from the
eigenvalues of $y_k$ by no more than $\sqrt{k} \cdot \epsilon$ in $\ell^2$
norm.  In particular this will happen if the $j$th terms of the sequences
differ by no more than $\epsilon.$ constructed in Lemma 2.1.  For
sufficiently large $k,$ $B_k \in \Gamma_R(x;m,k,\gamma).$

Suppose $\epsilon_0 > \epsilon >0$ and for each $k$ define $\Omega_k$ to
be the intersection of \[ [b_{1k} - \epsilon, b_{1k} + \epsilon ] \times
\cdots \times [b_{k k} - \epsilon, b_{k k} + \epsilon] \]

\noindent with $\mathbb R^k_<.$ Integrating according to the density given
above it follows that $\text{vol}(\Theta_{\epsilon}(B_k))$ exceeds
\begin{eqnarray} D_{k} \cdot \int_{\Omega_k } \Pi_{i <j}(t_i -t_j)^2 \,
dt_1 \cdots dt_{k}. \end{eqnarray} 

Recall the definitions of $S_k$ and $W_k$ in subsection 2.2 and denote by
$\Lambda_k$ the subset of $\Omega_k$ consisting of all $(t_1, \ldots,
t_k)$ satisfying $|t_i - t_j| \geq |b_{ik} - b_{jk}|.$ (3) dominates

\begin{eqnarray} D_{k} \cdot \Pi_{(i,j) \in
W_k} (b_{ik} - b_{jk})^2 \int_{\Lambda_k} \Pi_{(i,j) \in S_k} (t_i -
t_j)^2 \, dt_1 \cdots dt_k .  \end{eqnarray} 

Consider the map $F : [-\epsilon, \epsilon]^k \cap \mathbb R^k_<
\rightarrow \Lambda_k \subset \Omega_k$ given by $F(t_1,\ldots, t_k) =
(t_1 +b_{1k}, \ldots, t_k + b_{kk}).$ By a change of variables and
Selberg's integral formula 

\begin{eqnarray*} \epsilon^{k^2} \cdot \Pi_{j=1}^k \frac{\Gamma(j+1)  
\Gamma(j)^2}{\Gamma(k+j)} & < &\int_{[-\epsilon, \epsilon]^k} \Pi_{i<j}
(t_i - t_j)^2 \, dt_1 \cdots dt_k \\ & = & k! \cdot \int_{[-\epsilon,
\epsilon]^k \cap \mathbb R^k_<} \Pi_{i<j}(t_i - t_j)^2 \, dt_1 \cdots dt_k
\\ & < & k! \cdot (2 \epsilon)^{k^2 - 2 \cdot \#S_k -k}
\cdot\int_{[-\epsilon, \epsilon]^k \cup \mathbb R_<^k} \Pi_{(i,j) \in S_k}
(t_i - t_j)^2 \, dt_1 \cdots dt_k \\ & \leq & k! \cdot (2 \epsilon)^{k^2 -
2 \cdot \#S_k -k} \cdot\int_{\Lambda_k} \Pi_{(i,j) \in S_k} (t_i - t_j)^2
\, dt_1 \cdots dt_k. \\ \end{eqnarray*}

\noindent From this it follows that $\text{vol} \hspace{.03in}
(\Theta_{\epsilon}(B_k)) \geq (3) \geq (4) \geq C_k \cdot
\epsilon^{\#S_k}$ where 

\[C_k = D_{k} \cdot \Pi_{(i,j) \in W_k} (b_{ik} - b_{jk})^2 \cdot
(k!)^{-1} \cdot 2^{2 \cdot \#S_k +k - k^2} \cdot \Pi_{j=1}^k
\frac{\Gamma(j+1) \Gamma(j)^2}{\Gamma(k+j)} .\]

Thus for any $\epsilon_0 > \epsilon >0$ we have that
$P_{\epsilon}(\Theta(B_k)) \cdot L_k \cdot (3 \epsilon)^{k^2}\geq
\text{vol} \hspace{.03in} (\Theta_{\epsilon}(B_k)) > C_k \cdot \epsilon^{2
\cdot \#S_k +k}.$ For large enough $k$

\[ 2 \cdot \#S_k \leq 2 \left ( \frac{([c_1k] - \sqrt{k})^2}{2} +
\sum_{j=2}^s \frac{[c_jk]^2}{2} \right) \leq \left ( \sum_{j=1}^s [c_jk]^2
\right) - k \leq (1 - \alpha)k^2 -k\] 

\noindent whence, $2 \cdot \#S_k + k \leq (1 - \alpha)k^2.$ For any
$\epsilon_0 >\epsilon >0$ $P_{\epsilon}(\Theta(B_k)) \geq C_k \cdot
L_k^{-2} \cdot 3^{-k^2} \cdot \epsilon^{\alpha k^2}.$ Because
$\Theta(B_k)$ is locally isometric by Lemma 6.1 of [1] it follows that

\[ H^{\alpha k^2}_{\epsilon_0} (\Gamma_R(x;m,k,\gamma)) \geq
H^{\alpha k^2}_{\epsilon_0} (\Theta(B_k)) \geq C_k \cdot L_k^{-1} \cdot
3^{- k^2}.\] 

\noindent $m$ and $\gamma$ being arbitrary it follows from Lemma 4.4 of
[3], Lemma 4.1, and (1) from section 2.2 that 

\begin{eqnarray*} \mathbb H^{\alpha}(x) \geq
\mathbb H^{\alpha}_{\epsilon_0}(x) & = & - \log 3 + \liminf_{k \rightarrow
\infty} k^{-2} \cdot \left ( \log C_k - \log L_k \right) \\ & \geq & -
\log 3 - \frac{1}{2} \log 2 \pi e+ \liminf_{k \rightarrow \infty} \left
[k^{-2} \cdot\log C_k + \frac{1}{2} \cdot \log k \right ] \\ & \geq & -
\delta_0(x) \cdot \log 2 - \log 3 - \frac{1}{2} \log 2 \pi e + \int
\int_{\mathbb R^2 - D} \log |s-t| \, d\mu(s) d\mu(t) + \\ & & \liminf_{k
\rightarrow \infty} \left [ k^{-2} \cdot \log D_k \cdot \Pi_{j=1}^k
\frac{\Gamma(j+1) \Gamma(j)^2}{\Gamma(k+j)} + \frac{1}{2} \cdot \log k
\right ] \\ & \geq & - \delta_0(x)  \log 2 - \log 3 - \frac{1}{2} \log 2e
+ \frac{3}{4} + \int \int_{\mathbb R^2 - D} \log |s-t| \, d\mu(s) d\mu(t)  
+ \\ & & \liminf_{k \rightarrow \infty} k^{-2} \cdot \log \left
[\Pi_{j=1}^k \frac{\Gamma(j+1)  \Gamma(j)^2}{\Gamma(k+j)} \right] \\ & = &
- \delta_0(x) \log 2 - \frac{1}{2} \log 288e + \frac{3}{4} + \int
\int_{\mathbb R^2 -D} \log |s-t| \, d\mu(s) d\mu(t). \end{eqnarray*} \end{proof}

\begin{remark} Because $\mathbb H^1(x) = \chi(x) + \frac{1}{2} \log
\left(\frac{2}{\pi e}\right)$ it is clear that the lower bound of Lemma
2.1 is not sharp.  This is also clear from the reduction to
packing/covering number computations which in the microstate setting
introduce non-sharp estimates. \end{remark} 

\begin{example} Suppose $\mu = \sum_{j=1}^{\infty} \frac{1}{2^j} \cdot
\delta_{\frac{1}{j}}$ where $\delta_{r}$ is the the Dirac mass
concentrated at $r.$ If $x = x^* \in L^{\infty}\left (\langle
\frac{1}{j}\rangle_{j=1}^{\infty}, \mu \right )$ is the identity
multiplication operator then it follows from [1] and [4] that $\mathbb
H(x) = \delta_0(x) = 1 - \sum_{j=1}^{\infty} \frac{1}{4^j} = \frac{2}{3}.$
Moreover, $- \infty < \sum_{i \neq j} \left ( \frac{1}{2^{(i+j)}} \cdot
\log |\frac{1}{i} - \frac{1}{j}| \right ) < \infty$ so that by what
preceded, $- \infty < \mathbb H^{\frac{2}{3}}(x) < \infty.$ \end{example}

\section{Free Additivity} In this section suppose $x_1,\ldots,
x_n$ are selfadjoint elements of $M$ and that for each $1 \leq i \leq n$
$\mu_i$ is the Borel measure on $sp(x_i)$ induced by $\varphi.$ Set
$\alpha_i = \delta_0(x_i)$ and $\beta = \alpha_1 + \cdots + \alpha_n.$

\begin{lemma} If $\{x_1,\ldots,x_n\}$ is a freely independent family, then

\[ K_1 \leq H^{\beta}(x_1,\ldots,x_n) - \sum_{i=1}^n \int
\int_{\mathbb R^2 - D} \log |s-t| \, d\mu_i(s) d\mu_i(t) \leq K_2 \]

\noindent where $K_1 = -\frac{n}{2} \log 288 e + \frac{3n}{4} - \beta 
\log 2 $ and $K_2 = n \log 16 \sqrt{n} + \frac{n}{4}.$ \end{lemma}

\begin{proof} First for the lower bound on the difference.  
Suppose $m \in \mathbb N,$ $1 > \epsilon_0, \gamma >0,$ and $R$ exceeds
the maximum of the operator norms of any of the $x_i.$ By Corollary 2.14
of [5] there exists an $N \in \mathbb N$ such that if $k \geq N$ and
$\sigma$ is a Radon probability measure on $((M^{sa}_k(\mathbb C))_R)^n$
invariant under the $(U_k)^{(n-1)}$-action given by $(\xi_1,\ldots, \xi_n)  
\mapsto (\xi_1, u_1 \xi_2 u_1^*, \ldots, u_{n-1}\xi_n u_{n-1}^*)$ where
$(u_1,\ldots, u_{n-1}) \in (U_k)^{(n-1)},$ then $\sigma(\omega_k) >
\frac{1}{2}$ where conclusion. \[ \omega_k =
\{(\xi_1,\ldots,\xi_n)\in((M_k^{sa}(\mathbb C))_{R+1})^{n} : \langle
\{\xi_i\} \rangle_{i=1}^n \hspace{.05in} \text{are} \left (m,
\frac{\gamma}{4^m} \right) \text{ - free} \}. \]

\noindent The preceding section provided for each $i$ a sequence $\langle B_{ik} \rangle_{k=1}^{\infty}$ such that
for any $m^{\prime} \in \mathbb N$ and $\gamma^{\prime} >0$ $B_{ik} \in \Gamma_R(x_i;m^{\prime},k,\gamma^{\prime})$
for sufficiently large $k.$ Also for any $k,$ $\|B_{ik}\| \leq R.$ Write $\Theta (B_{ik})$ for the unitary orbit of
$B_{ik}$ and $g_{ik}$ for the topological dimension of this orbit.  The proof of Lemma 2.1 yielded constants
$C_{ik}$ for each $1 \leq i \leq n$ such that for any $\epsilon >0$

\[ P_{\epsilon}(\Theta(B_{ik})) \geq C_k \cdot L_k^{-1} \cdot 3^{-k^2} 
\epsilon^{\#S_{ik} 
- k^2} 
\]

\noindent where $S_{ik} \leq (1 - \alpha_i)k^2.$ For each $k \in \mathbb
N$ denote by $\mu_k$ the probability measure of $((M^{sa}_k(\mathbb
C)_{R+1})^n$ obtained by restricting $\sum_{i=1}^n g_{ik}$-Hausdorff
measure (with respect to the $|\cdot|_2$ norm) to the $\sum_{i=1}^n
g_{ik}$ -dimensional manifold $T_k = \Theta(B_{1k}) \times \cdots \times
\Theta (B_{nk})$ and normalizing appropriately.  $\mu_k$ is a Radon
probability measure invariant under the $(U_k)^{n-1}$-action described
above because such actions are isometric, whence $\mu_k(\omega_k) >
\frac{1}{2}.$ Set $F_k = \omega_k \cap T_k.$ It is clear that $\mu_k(F_k)
= \mu_k(\omega_k) > \frac{1}{2}$ and for large enough $k,$ $F_k \subset
\Gamma_{R+1}(x_1,\ldots, x_n;m,k,\gamma).$

$T_k$ is a locally isometric, smooth, compact manifold of dimension 
$\sum_{i=1}^n g_{ik}$ (by locally isometric we means that for any 
$\epsilon >0$ any two open $\epsilon$ balls of the metric space are 
isometric).  Moreover by the preceding paragraph 

\[ P_{\epsilon}(T_k) \geq \Pi_{i=1}^n P_{\epsilon}(\Theta(B_{ik})) \geq 
\Pi_{i=1}^n  C_{ik} \cdot L_k^{-1} \cdot 3^{-k^2} \epsilon^{\#S_{ik} - 
k^2}. \]

\noindent This estimate holds for all $1 > \epsilon >0.$  By Lemma 6.1 of [1] 

\[ H^{\beta k^2}_{\epsilon_0}(\Gamma_R(x_1,\ldots, x_n;m,k,\gamma)) \geq
H^{\beta k^2}_{\epsilon_0} (F_k) \geq H^{\sum_{i=1}^n (k^2 -
S_{ik})}_{\epsilon_0}(F_k) \geq L_k^{-n} \cdot 3^{-nk^2} \cdot \frac{1}{2}
\cdot \Pi_{i=1}^n C_{ik}.\]

\noindent Thus using the computations already made in Lemma 3.2,

\begin{eqnarray*} \mathbb H^{\beta}_{\epsilon_0}(x_1,\ldots, x_n;m,\gamma)
& \geq & \liminf_{k \rightarrow \infty} k^{-2} \log \left (L_k^{-n} \cdot
3^{-nk^2} \cdot \frac{1}{2} \cdot \Pi_{i=1}^n C_{ik} \right) 
\\ & \geq &\sum_{i=1}^n \left[ - \log 3 + \liminf_{k \rightarrow \infty} k^{-2} \cdot
\log (C_{ik} - \log L_k) \right] \\ & \geq & -\frac{n}{2} \log 288e +
\frac{3n}{4} - \beta \log 2 + \sum_{i=1}^n \int \int_{\mathbb R^2 -D} \log 
|s-t|
d\mu_i(s)  d\mu_i(t).\\ \end{eqnarray*}

Now for the upper bound.  Suppose $\epsilon, t >0.$ Recalling the proof
of Lemma 3.1 we can produce $m \in \mathbb N$ and $\gamma >0$ such that

\begin{eqnarray*} H^{\beta k^2}_{\epsilon}(\Gamma(x_1,\ldots,
x_n;m,k,\gamma)) & \leq & H^{\beta k^2}_{\epsilon} \left (\Pi_{i=1}^n
\Theta_t(A_{ik}) \right) \\ & \leq & K_{\frac{\epsilon}{2}} \left
(\Pi_{i=1}^n \Theta_t(A_{ik}) \right) \cdot \epsilon^{\beta k^2} \\ &
\leq & \left [ \Pi_{i=1}^n
K_{\frac{\epsilon}{2\sqrt{n}}}(\Theta_t(A_{ik})) \right] \cdot
\epsilon^{\beta k^2} \\ & \leq &
\left [ \Pi_{i=1}^n P_{\frac{\epsilon}{4\sqrt{n}}} (\Theta_t(A_{ik}))  
\right] \cdot \epsilon^{\beta k^2} \\ & \leq & 
\left [ \Pi_{i=1}^n \text{vol} (\Theta_{t + \frac{\epsilon}{4 
\sqrt{n}}}(A_{ik})) \right] \cdot \frac{(4 \sqrt{n})^{n k^2} 
\epsilon^{(\beta -n)k^2}}{L_k^n}.
\end{eqnarray*}

\noindent Thus,

\[ \mathbb H^{\beta}_{\epsilon}(x_1,\ldots, x_n;m,\gamma) \leq
\sum_{i=1}^n \limsup_{k \rightarrow \infty} k^{-2} \cdot
[\log(\text{vol}\hspace{.03in} \Theta_{t + \frac{\epsilon}{4
\sqrt{n}}}(A_{ik})) - \log L_k] + \log (4 \sqrt{n}) + (\alpha_i -1 )
\log \epsilon. \]

\noindent This bound being independent of $m$ and $\gamma,$ letting
$\epsilon \rightarrow 0,$ and using the computation already made in Lemma
3.1, we conclude that

\[ \mathbb H^{\beta}(x_1,\ldots, x_n) \leq \left ( \sum_{i =1}^n \int
\int_{\mathbb R^2 - D} \log |s-t| \,d\mu_i(s) d\mu_i(t) \right )+ n \log
16 \sqrt{n} + \frac{n}{4}.\]

\end{proof}


\begin{thebibliography}{[ASMR]}

\bibitem{1} Jung, Kenley {\it Fractal dimensions and entropies for
microstate spaces}, preprint.

\bibitem{2} Mehta, M.L. {\it Random Matrices}, Academic Press, 1991.

\bibitem{3} Voiculescu, D. {\it The analogues of entropy and of Fisher's
information measure in free probability theory, II}, Inventiones
mathematicae 118, (1994), 411-440.

\bibitem{4} Voiculescu, D. {\it The analogues of entropy and of Fisher's
information measure in free probability theory III: The Absence of Cartan
Subalgebras}, Geometric and Functional Analysis, Vol. 6, No.1 (1996)
(172-199).

\bibitem{5} Voiculescu, D. {\it A strengthened
asymptotic freeness result for random matrices with applications to free
entropy}, IMRN, 1 (1998), 41-64.

\end{thebibliography}
\end{document}